\documentclass{iopart}

\usepackage{iopams}
\usepackage{stackrel}
\usepackage[utf8]{inputenc}
\usepackage{amssymb}
\usepackage{amsthm}
\usepackage{tikz}
\usepackage{mathrsfs}
\usepackage{epsf}
\usepackage{bookmark,hyperref}
\hypersetup{
     colorlinks=true,
     linkcolor=blue,
     filecolor=blue,
     citecolor = blue,
     urlcolor=cyan,
     }
\usepackage{pstricks}

\newcommand{\R}{{\mathbb R}}

\newcommand{\I}{{\mathcal{O}}}

\newcommand{\dem}{{\em Proof: \;}}
\newcommand{\fdem}{\hfill $\square$}

\newtheorem{teo}{Theorem}[section]
\newtheorem{lema}[teo]{Lemma}
\newtheorem{cor}[teo]{Corollary}

\newtheorem*{teosn}{Theorem}

\def\set#1{\left\{\, #1 \,\right\}}

\begin{document}

\title[Lyapunov instability in Newtonian dynamics]{On the Lyapunov instability in Newtonian dynamics}

\author{J. M. Burgos$^1$, E. Maderna$^2$ and M. Paternain$^3$}

\address{$^1$ Departamento de Matem\'aticas, CINVESTAV--CONACYT, Av. Instituto Polit\'ecnico Nacional 2508, Col. San Pedro Zacatenco, 07360 Ciudad de M\'exico, M\'exico.}
\address{$^2$ IMERL, Facultad de Ingenier\'ia, Universidad de la Rep\'ublica, 
Av. Herrera y Reissig 565, 11300 Montevideo, Uruguay.}
\address{$^3$ Centro de Matem\'atica, Facultad de Ciencias, Universidad de la Rep\'ublica, Igu\'a 4225,
11400 Montevideo, Uruguay.}
\eads{\mailto{burgos@math.cinvestav.mx}, \mailto{eze@fing.edu.uy}, \mailto{miguel@cmat.edu.uy}}

\begin{abstract}
We prove Lyapunov instability for cases in which the local minimum of the potential energy is reached on a hypersurface of the configuration space.
In contrast to the known results in this direction, which hold for potentials satisfying hypotheses on the first non-zero jet, this new result covers several real analytic cases that the previous do not.
\end{abstract}

\maketitle

\section{Introduction}

The Lagrange-Dirichlet Theorem concerns conservative holonomic mechanical systems with finite degrees of freedom and states that every strict local minimum of the potential is a Lyapunov stable equilibrium point of the dynamics. It was stated by Lagrange in \cite{Lagrange} and proved by Dirichlet in \cite{Dirichlet}.

However, without further hypothesis besides differentiablity, the converse is false. In 1904 Painlev\'e proposed the following counterexample (see \cite{Pai}): Consider the one degree of freedom mechanical system
$\ddot{x}= -\nabla U(x)$ with the $C^{\infty}$ potential
$$U(x)= \exp(-|x|^{-1})\sin(|x|^{-1})\,.$$
The origin is a critical point and it is not a minimum. For every neighborhood of the origin, there is an interval centered at it and contained in the neighborhood such that the potential is maximum and strictly greater than zero on the interval's boundary. Therefore, for every neighborhood of the origin, every motion with small enough energy is trapped in an interval contained in the neighborhood hence the origin is Lyapunov stable and a counterexample of the Lagrange-Dirichlet converse.

A more striking example is the following proposed by Laloy in \cite{Laloy1}. Consider the two degrees of freedom mechanical system $\ddot{{\bf x}}= -\nabla U({\bf x})$, ${\bf x}=(x,y)$, with the $C^{\infty}$ potential
$$U(x,y)= \exp(-|x|^{-1})\sin(|x|^{-1})- \exp(-|y|^{-1})\sin(|y|^{-1}) - y^{2}\,.$$
The origin is a critical point and it is not a minimum. In contrast with the previous example, now there are no trapping zones for the set $U^{-1}\left((-\infty, 0)\right)$ contains the two diagonals $x=\pm y$ and these are the only escape routes to infinity where a priory any motion starting near the origin could take. However, the projection of the motion on the first coordinate is governed by the previous example hence the origin is again Lyapunov stable and another counterexample of the Lagrange-Dirichlet converse.

Then, a natural question arises: What conditions are needed in order for the Lagrange-Dirichlet converse to hold? In this respect, we find the following in Arnold's Problems book \cite{ArnoldProblem}:
\begin{quote}
``\textbf{1971-4.} Prove the instability of the equilibrium $\textbf{0}$ of the analytic system $\ddot{x}=-\partial U/\partial x$ in the case where the isolated critical point $\textbf{0}$ of the potential $U$ is not a minimum."
\end{quote}

For the bibliography and comments related to the problem we highly recommend the comment section in \cite{ArnoldProblem}, pages 250-253. In that section, the authors recognize the fact that:
\begin{quote}
``...The problem on the converse of the Lagrange-Dirichlet theorem makes therefore sense only under one or another additional assumptions (e.g., that of analyticity of the potential)."
\end{quote}
Lyapunov himself stated the problem for real analytic potentials in \cite{Ly}.

In \cite{Br}, Brunella solves the Arnold's problem for two degrees of freedom (Corollary in \cite{Br}, page 1346.).

In \cite{Palamodov}, Palamodov completely solves the Arnold's problem giving a beautiful proof in terms of real analytic geometry using monoidal transformations also known as \textit{blow-ups}. Concretely, he proves (Corollary 2.2 in \cite{Palamodov}, page 7.): Let $U$ be a real analytic potential and $p$ a critical point. If $p$ belongs to the closure of the region where the potential is strictly less than $U(p)$, then $p$ is Lyapunov unstable.

Palamodov proves the Theorem in the context of Lagrangian dynamics with a mechanical Lagrangian (or natural system as he calls it). In particular, to prove the Lagrange-Dirichlet converse, only the case of a non strict minimum critical point is left.

Starting from Lyapunov \cite{Ly} and following many others \cite{Ref1}, \cite{Ha}, \cite{Ref2}, \cite{Ref3}, \cite{Ref4}, \cite{Ref5}, \cite{Ref6}, \cite{Ref7}, many partial results have been given towards this direction and their common thing is that the Lyapunov instability criteria involves the lack of a local minimum at the origin of the first nonzero $k$--th order jet of the potential with $k\geq 2$ \footnote{In \cite{Ref3}, the degenerate second order differential is admitted.}. However, these criteria are not sufficient to prove the case of a non strict local minimum of the potential. As an example, consider the \textit{gutter} potential $U(x,y)= x^{4}$ and see that none of the instability criteria described before apply. However, any non trivial motion escapes through $x=0$ hence any critical point is Lyapunov unstable.

The case of a non positive potential is trivial because in this case for every positive energy the corresponding Jacobi--Mapertuis metric is complete and by the Hopf--Rinow Theorem there is a trajectory from the critical point to any other point with arbitrarily small energy.

The first open problem described in section 3, \textit{Open problems and a conjecture}, in \cite{Palamodov} is the study of a non strict local minimum of a real analytic potential. As far as we know, it is still open. The case of a non strict local minimum of the potential but with two degrees of freedom was treated in \cite{Laloy}.

In this note we restrict ourselves to Newtonian dynamics, i.e. $\ddot{x}= -\nabla U(x)$, and study the case of a non strict local minimum hypersurface of the following class of potentials for an arbitrary number of degrees of freedom:

\bigskip
\noindent\textbf{Hypothesis.}
{\it The potential $U$ is the composition $g\circ f$ such that zero is a regular value of $f$ in $C^{3}(\R^{n}, \R)$ and $g$ in $C^{2}(\R, \R)$ verifies $g\geq 0$ vanishing only at zero.}
\vspace{4pt}

In particular the set 
$M=\set{x\in\R^n\mid f(x)=0}=\set{x\in\R^n\mid U(x)=0}$
is a hypersurface of $\R^{n}$ and $M\times\{\bf 0\}$ consists entirely of equilibrium points of the Newtonian dynamics.

\begin{teosn}\label{Main}
Every point in $M\times\{\bf 0\}$ is Lyapunov unstable.
\end{teosn}

It is worth to mention that our result covers several potentials that the previous do not. As an example, consider the following potential:
$$U(x,y,z)= \left(x^{2}+2 y^{2}+ 3 z^{2}-1\right)^{4}.$$
Its zero potential critical locus is an ellipsoid and it is a minimum of the potential. Its first nonzero $k$--th order jet with $k\geq 2$ at every critical point is with $k=4$ and has a local minimum at the origin hence none of the mentioned analytic methods can be applied.

\section{Proof of the instability}

Let $p$ be a point in $M$, $v$ a non zero vector in $T_p M$ and for every $\varepsilon>0$ consider the solution $x^{\varepsilon}$ of the Cauchy problem
\begin{equation}\label{Cauchy1}
\ddot{x}^{\varepsilon}= -\nabla U(x^{\varepsilon}),\ x^{\varepsilon}(0)=p,\ \dot{x}^{\varepsilon}(0)= \varepsilon v,\ \varepsilon>0.
\end{equation}

For every $\varepsilon>0$, define $x_\varepsilon$ such that $x_\varepsilon(\tau)= x^{\varepsilon}(\tau/\varepsilon)$ where $x^{\varepsilon}$ defined. These are solutions of the  Cauchy problem:
\begin{equation}\label{Cauchy2}
\ddot{x}_{\varepsilon}= -\varepsilon^{-2}\nabla U(x_{\varepsilon}),\ x_{\varepsilon}(0)=p,\ \dot{x}_{\varepsilon}(0)= v,\ \varepsilon>0.
\end{equation}
Now, the initial conditions are fixed but the equation becomes singular as $\varepsilon\to 0^{+}$. Denote by $I_\varepsilon$ the maximal interval containing zero where $x_{\varepsilon}$ is defined.

\begin{lema}\label{Lema1}
For every $\varepsilon>0$, $\Vert\dot{x}_\varepsilon(\tau)\Vert\leq \Vert v \Vert$ for every $tau$ in $I_\varepsilon$ and
$${\rm Im}(x_\varepsilon)\subset [U\leq \varepsilon^{2}\Vert v \Vert^{2}/2].$$
\end{lema}
\dem
For every $\varepsilon>0$, the Hamiltonian
$$H_\varepsilon({\rm x},{\rm v})= \Vert {\rm v} \Vert^{2}/2+\varepsilon^{-2}U({\rm x})$$
is constant along the solution $x_\varepsilon$ hence
$$\Vert \dot{x}_\varepsilon(\tau) \Vert^{2}/2,\ \varepsilon^{-2}U(x_\varepsilon(\tau))\leq H_\varepsilon(x_\varepsilon(\tau), \dot{x}_\varepsilon(\tau))= H_\varepsilon(p, v)= \Vert v \Vert^{2}/2.$$
\fdem

\begin{cor}\label{Cor1}
Let $T>0$. For every $\varepsilon>0$ and every $\tau$ in $I_\varepsilon\cap [-T, T]$, 
$$(x_{\varepsilon}(\tau), \dot{x}_\varepsilon(\tau))\in \overline{B(p, T \Vert v \Vert)}\times \overline{B({\bf 0}, \Vert v \Vert)}.$$
Note that the region is a compact set not depending on $\varepsilon$.
\end{cor}
\dem
By Lemma \ref{Lema1}, $\Vert \dot{x}_\varepsilon(\tau)\Vert\leq \Vert v \Vert$ and
$$\Vert x_\varepsilon(\tau)-p \Vert= \left\Vert\int_0^{\tau}ds\ \dot{x}_\varepsilon(s)\right\Vert \leq \left\vert\int_0^{\tau}ds\ \Vert\dot{x}_\varepsilon(s)\Vert\ \right\vert\leq |\tau|\ \Vert v \Vert\leq T \Vert v \Vert,$$
the result follows.
\fdem

\begin{cor}\label{Cor2}
For every $\varepsilon>0$, $x_\varepsilon$ is defined over the whole real line.
\end{cor}
\dem
Consider the maximal interval $I_\varepsilon=(\omega_-, \omega_+)$ and suppose that $\omega_+$ is finite. Then, $(x_\varepsilon, \dot{x}_\varepsilon)|_{[0,\omega_+)}$ is contained in the compact set
$$\overline{B(p, \omega_+ \Vert v \Vert)}\times \overline{B({\bf 0}, \Vert v \Vert)}$$
which is absurd hence $\omega_+=+\infty$. Analogously, $\omega_-=-\infty$.
\fdem

\begin{cor}\label{Cor3}
Let $T>0$. There is a continuous curve $x:[-T,T]\rightarrow M$ with $x(0)=p$ and a sequence $(\varepsilon_j)$ such that $\varepsilon_j>0$, $\varepsilon_j\to 0^{+}$ and $x_{\varepsilon_j}\rightarrow x$ uniformly on $[-T,T]$.
\end{cor}
\dem
Because $x_\varepsilon(0)=p$ for every $\varepsilon>0$, by Arzel\`a--Ascoli Theorem, there is such a sequence and a continuous curve $x:[-T,T]\rightarrow \R^{n}$ such that $x_{\varepsilon_j}\rightarrow x$ uniformly on $[-T,T]$. For every $j$, ${\rm Im}(x_{\varepsilon_j})$ is contained in $\left[U\leq \varepsilon_j^{2}\Vert v \Vert^{2}/2\right]$ and this is a nested sequence. Then,
$${\rm Im}(x)\subset \bigcap_{j}\left[U\leq \varepsilon_j^{2}\Vert v \Vert^{2}/2\right]= [U=0]=M.$$
\fdem

Now we construct suitable coordinates in order to isolate the singular limit in \eref{Cauchy2} in one coordinate. Consider the flow $\phi$ in $\R^{n}- Crit(f)$
\begin{equation}\label{Cauchy3}
\partial_t \phi= \frac{\nabla f}{\Vert \nabla f \Vert^{2}}(\phi),\ \phi(0,x)= x,\ x\in \R^{n}-Crit(f).
\end{equation}

Consider a local coordinate neighborhood $(V, \psi)$ of $M$ centered at $p$ and denote $w$ the velocity in these coordinates:
\begin{equation}\label{Coord_w}
d\psi({\bf 0},w)= (p,v).
\end{equation}
By Hypothesis, $M$ is contained in $\R^{n}-Crit(f)$. We define the $C^{2}$ local coordinate neighborhood $(\I(V), \Psi)$ such that
$$\Psi(r,y)= \phi(r, \psi(y))$$
where $\I(V)\subset \R^{n}-Crit(f)$ is the union of the set of orbits of \eref{Cauchy3} with initial condition in $V$.

\begin{lema}\label{Lema2}
\begin{enumerate}
\item $\Psi(0,y)= \psi(y)$ for every $y$ in $\psi^{-1}(V)$.
\item $U(\Psi(r,y))= g(r)$ for every $(r,y)$ in $\Psi^{-1}(\I(V))$.
\item $(\I(V), \Psi)$ is a local coordinate neighborhood of $\R^{n}$.
\end{enumerate}
\end{lema}
\dem
\begin{enumerate}
\item $\Psi(0,y)= \phi(0, \psi(y))= \psi(y)$.
\item By definition, $\partial_t(f\circ\phi)= \langle\nabla f(\phi), \partial_t\phi\rangle= 1$ hence
$$f\circ\phi(t, x)= t+ f(\phi(0,x))= t+f(x),$$
$$f(\Psi(r,y))= r+f(\psi(y))=r$$
for $\psi(y)$ is in $M=[f=0]$. Then, $U(\Psi(r,y))= g\circ f(\Psi(r,y))= g(r)$.

\item Define $\phi_t(x):=\phi(t,x)$. Because $\nabla f(\psi(y))\neq {\bf 0}$ and $\nabla f(\psi(y))\perp T_{\psi(y)}V$, $d_{(0,y)}\Psi$ is an isomorphism  for $d_{y}\psi$ is so. By definition of $\Psi$, we have the relation $\Psi(r,y)= \phi_r\circ\Psi(0,y)$. Taking differentials,
$$d_{(r,y)}\Psi= d_{\psi(y)}\phi_r\circ d_{(0,y)}\Psi.$$
Then, $d_{(r,y)}\Psi$ is also an isomorphism for $d_{\psi(y)}\phi_r$ is so by the Liouville formula. By the inverse function Theorem, $\Psi$ is a local diffeomorphism. To show that $\Psi$ is an embedding, it rest to show that it is injective.

Suppose that $\Psi(r,y)=\Psi(r', y')$. Then,
$$r= f(\Psi(r,y))= f(\Psi(r',y'))= r',$$
$$\psi(y)= \phi_{-r}(\Psi(r,y))= \phi_{-r}(\Psi(r,y'))= \psi(y')$$
so $(r,y)=(r', y')$ for $\psi$ is injective.
\end{enumerate}
\fdem

Let $T>0$ be small enough such that the compact set $\overline{B(p, T\Vert v \Vert)}$ is contained in $\I(V)$. From now on, all the curves will be defined on $[-T,T]$.

The chart $(\I(V), \Psi)$ defines $C^{2}$ curvilinear coordinates with versors $e_r,\ e_k$ and scale factors $h_r, h_k>0$ respectively:
$$\partial_r \Psi= h_r\ e_r,\ \partial_k \Psi= h_k\ e_k.$$
Denote by $e^{r}, e^{k}$ the duals of $e_r,\ e_k$ respectively\footnote{Because $e_r\perp e_k$, $k=1, \ldots, n-1$, we have $e_r=e^{r}$.}. Denote by $r_\varepsilon$ and $y_\varepsilon$ the coordinates of $x_\varepsilon$ with respect to the chart $(\I(V), \Psi)$:
$$\Psi(r_\varepsilon(\tau), y_\varepsilon(\tau))= x_\varepsilon(\tau)$$
and an analogous definition for the limit curve in Corollary \ref{Cor3}.

\begin{lema}\label{Lema3}
\begin{enumerate}
\item $r_\varepsilon\rightarrow 0$ uniformly on $[-T,T]$ as $\varepsilon\to 0^{+}$.
\item $\dot{r}_\varepsilon, \dot{y}_\varepsilon^{k}$ are bounded by a constant not depending on $\varepsilon$.
\end{enumerate}
\end{lema}
\dem
\begin{enumerate}
\item $g(r_\varepsilon)= U(x_\varepsilon)\leq \varepsilon^{2}\Vert v \Vert^{2}/2\rightarrow 0$ as $\varepsilon\to 0^{+}$.
\item For every ${\rm x}$ in $\Psi^{-1}(\I(V))$ define the quadratic form $Q_{{\rm x}}$ such that
$$Q_{{\rm x}}({\rm v})= \Vert d_{{\rm x}} \Psi({\rm v})\Vert^{2}.$$
It is positive definite for every ${\rm x}$ and defines a strictly positive continuous function on the unit tangent sphere bundle $\pi:T^{1}\R^{n}\to \R^{n}$. In particular, it attains a minimum value $m>0$ on the compact set $(\Psi\circ\pi)^{-1}(\overline{B(p,T\Vert v \Vert )})$.
For every $\varepsilon>0$ and $\tau$ in $[-T,T]$ we have
$$m\Vert(\dot{r}_{\varepsilon}(\tau), \dot{y}_{\varepsilon}(\tau))\Vert^{2}\leq Q_{(r_{\varepsilon}(\tau), y_{\varepsilon}(\tau))}(\dot{r}_{\varepsilon}(\tau), \dot{y}_{\varepsilon}(\tau))= \Vert \dot{x}_\varepsilon(\tau)\Vert^{2}\leq \Vert v \Vert^{2}$$
hence
$$|\dot{r}_{\varepsilon}(\tau)|,\ |\dot{y}_{\varepsilon}^{k}(\tau))|\leq m^{-1/2}\Vert v \Vert.$$
and the result follows.
\end{enumerate}
\fdem

\begin{lema}\label{Lema4}
For every $\varepsilon>0$, considering the functions $r_\varepsilon$ and $\dot{r}_\varepsilon$ as external parameters, we have the nonautonomous equations
\begin{equation}\label{Equation_y}
\fl\frac{h_r}{h_k} \bigg\langle e^{k}, \frac{\partial e_r}{\partial r}\bigg\rangle \dot{r}_\varepsilon^{2} + \frac{1}{h_k}\bigg\langle e^{k}, \frac{\partial^{2}\Psi}{\partial y^{a} \partial y^{b}}\bigg\rangle \dot{y}_\varepsilon^{a}\dot{y}_\varepsilon^{b} +2\frac{h_r}{h_k} \bigg\langle e^{k}, \frac{\partial e_r}{\partial y^{a}}\bigg\rangle \dot{y}_\varepsilon^{a}\dot{r}_\varepsilon +\ddot{y}_\varepsilon^{k}=0
\end{equation}
where the coefficients are evaluated over $y_\varepsilon$ and $r_\varepsilon$.
\end{lema} 
\dem
For every $\tau$ in $[-T, T]$, because of our hypothesis, either $\nabla U(x_\varepsilon(\tau))={\bf 0}$ or $\nabla U(x_\varepsilon(\tau))$ is collinear with $\nabla f(x_\varepsilon(\tau))$. On the other hand, $\nabla f$ is collinear with $e_r$ at every point in $\R^{n}-Crit(f)$. Thus, with respect to the $(\I(V), \Psi)$ coordinate chart, the motion equation \eref{Cauchy2} reads as follows
\begin{equation}\label{SecondNewtonLaw}
\ddot{x}_\varepsilon(s)= -\varepsilon^{-2}\Vert \nabla U(x_\varepsilon)\Vert\ e_r
\end{equation}
where the acceleration has the following expression

\begin{equation}\label{AccelerationCurvilinear}
\ddot{x}_\varepsilon= \frac{\partial^{2}\Psi}{\partial r^{2}}\dot{r}_\varepsilon^{2} + \frac{\partial^{2}\Psi}{\partial y^{a}\partial y^{b}}\dot{y}_\varepsilon^{a}\dot{y}_\varepsilon^{b}+2\frac{\partial^{2}\Psi}{\partial y^{a}\partial r}\dot{y}_\varepsilon^{a}\dot{r}_\varepsilon + \frac{\partial\Psi}{\partial r}\ddot{r}_\varepsilon+\frac{\partial\Psi}{\partial y^{a}}\ddot{y}_\varepsilon^{a}.
\end{equation}

\noindent In terms of the scale factors and versors, expression \eref{AccelerationCurvilinear} becomes
\begin{eqnarray*}
\fl\ddot{x}_\varepsilon = \left(\frac{\partial h_r}{\partial r}\dot{r}_\varepsilon^{2} + h_r \bigg\langle e^{r}, \frac{\partial e_r}{\partial r}\bigg\rangle \dot{r}_\varepsilon^{2} + \bigg\langle e^{r}, \frac{\partial^{2}\Psi}{\partial y^{a}\partial y^{b}}\bigg\rangle\dot{y}_\varepsilon^{a}\dot{y}_\varepsilon^{b}\right. \\
\left.+2h_r \bigg\langle e^{r}, \frac{\partial e_r}{\partial y^{a}}\bigg\rangle \dot{y}_\varepsilon^{a}\dot{r}_\varepsilon+ 2\frac{\partial h_r}{\partial y^{a}}\dot{y}_\varepsilon^{a}\dot{r}_\varepsilon + h_r \ddot{r}_\varepsilon\right)e_r \\
\fl+\left(h_r \bigg\langle e^{k}, \frac{\partial e_r}{\partial r}\bigg\rangle \dot{r}_\varepsilon^{2} + \bigg\langle e^{k}, \frac{\partial^{2}\Psi}{\partial y^{a}\partial y^{b}}\bigg\rangle \dot{y}_\varepsilon^{a}\dot{y}_\varepsilon^{b} +2h_r \bigg\langle e^{k}, \frac{\partial e_r}{\partial y^{a}}\bigg\rangle \dot{y}_\varepsilon^{a}\dot{r}_\varepsilon +h_k \ddot{y}_\varepsilon^{k}\right)e_k.
\end{eqnarray*}
Combining the motion equation \eref{SecondNewtonLaw} with expression \eref{AccelerationCurvilinear} we obtain the equations \eref{Equation_y}.
\fdem

Now, the equations of motion \eref{Equation_y} are not singular as $\varepsilon\to 0^{+}$.

\begin{cor}\label{Cor4}
There is a constant $C$ not depending on $\varepsilon$ such that $\Vert\ddot{y}_\varepsilon(\tau)\Vert\leq C$ for every $\varepsilon>0$ and every $\tau$ in $[-T,T]$.
\end{cor}
\dem
All of the coefficients are continuous on $\Psi^{-1}(\I(V))$ hence they are bounded on the compact set $\Psi^{-1}(\overline{B(p,T\Vert v \Vert )})$. By Lemma \ref{Lema3}, all of the velocities are bounded by a constant not depending on $\varepsilon$ therefore, by Lemma \ref{Lema4}, the same occurs with the accelerations.
\fdem

\begin{cor}\label{Cor5}
Taking a subsequence if necessary of the sequence in Corollary \ref{Cor3}, the curve $x$ is in $C^{1}[-T,T]$ with $\dot{x}(0)=v$.
\end{cor}
\dem
Recall that the coordinates are centered at $p$ and $w$ is the initial velocity with respect to these, see the equation \eref{Coord_w}. Because $\dot{y}_\varepsilon(0)=w$ for every $\varepsilon>0$, by the previous Corollary and Arzel\`a--Ascoli Theorem, taking a subsequence if necessary of the sequence in Corollary \ref{Cor3}, we have $\dot{y}_{\varepsilon_j}\rightarrow e$ uniformly on $[-T,T]$ for some continuous function $e$ such that $e(0)=w$. For every $j$ we have
$$y_{\varepsilon_j}(\tau)= \int_0^{\tau} ds\ \dot{y}_{\varepsilon_j}(s)$$
and taking the limit as $j\to +\infty$,
$$y(\tau)= \int_0^{\tau} ds\ e(s).$$
In particular, $\dot{y}=e$ is continuous and $\dot{y}(0)=w$. Because $r_{\varepsilon_j}\to 0$ as $j\to +\infty$ and $x= \Psi(0,y)=\psi(y)$, we have the result.
\fdem

\bigskip

\noindent\textit{Proof of Theorem \ref{Main}:}
By Corollaries \ref{Cor3} and \ref{Cor5}, there is a $C^{1}$ curve $x:[-T,T]\rightarrow M$ with $T>0$, $x(0)=p$, $\dot{x}(0)=v$ and a sequence $(\varepsilon_j)$ such that $\varepsilon_j>0$, $\varepsilon_j\to 0^{+}$ and $x_{\varepsilon_j}\rightarrow x$ uniformly on $[-T,T]$.

\noindent Consider the continuous function $(\tau\mapsto \Vert x(\tau)-p \Vert)$ on $[0,T]$. Because $x$ is differentiable at $\tau=0$ with $\dot{x}(0)= v\neq {\bf 0}$, the function attains a maximum $R>0$ at $\tau^{*}$ in $(0,T]$.

\noindent There is a natural $j_0$ such that $\Vert x_{\varepsilon_j}(\tau^{*})-x(\tau^{*})\Vert <R/2$ if $j\geq j_0$. In particular,
$$x^{\varepsilon_j}(\tau^{*}/\varepsilon_j)\notin B(p, R/2),\ \ \ j\geq j_0$$
while $(x^{\varepsilon_j}(0), \dot{x}^{\varepsilon_j}(0))= (p, \varepsilon_j v)\rightarrow (p,{\bf 0})$ as $j\to +\infty$. We conclude that $(p,{\bf 0})$ is a Lyapunov unstable equilibrium point.

\noindent Because the choice of the point $p$ was arbitrary, we have the result.
\fdem

\ack
We are very grateful with the anonymous referees for their work, their suggestions considerably improved the paper. The first author has a CONACYT research fellowship. The third author was supported by the FCE-ANII-135352 grant.

\References

\bibitem[Ar]{ArnoldProblem}
Arnold V I 2002 \emph{Arnold's Problems}, Springer-Verlag.

\bibitem[Br]{Br}
Brunella M 1998 \emph{Instability of equilibria in dimension three}, Ann I Fourier {\bf 48}.

\bibitem[Di]{Dirichlet}
Dirichlet L G 1846 \emph{\"Uber die Stabilitat des Gleichgewichts}, J. Reine Angew. Math. {\bf 32} 85--88.

\bibitem[GT]{Ref1}
Garcia M V P, Tal F A 2003 \emph{Stability of equilibrium of conservative systems with two degrees of freedom}, J. Differential Equations {\bf 194} 364--81.

\bibitem[Ha]{Ha}
Hagedorn P 1971 \emph{Die Umkehrung der Stabilit\"atss\"atze von Lagrange-Dirichlet und Routh}, Arch. Rational Mech. Anal. {\bf 42} 281--316.

\bibitem[Ko1]{Ref2}
Kozlov V V 1982 \emph{Asymptotic solutions of equations of classical mechanics}, J. Appl. Math. Mech. \textbf{46} 454--7.

\bibitem[Ko2]{Ref3}
------- 1987 \emph{Asymptotic motions and the inversion of the Lagrange-Dirichlet theorem}, J. Appl. Math. Mech. \textbf{50} 719--25.

\bibitem[Ko3]{Pai}
------- 1995 \emph{Problemata nova, ad quorum solutionem mathematici invitantur}, Transl. Amer. Math.Soc. Ser. 2 {\bf 168} 141--171.

\bibitem[Ku]{Ref4}
Kuznetsov A N 1989 \emph{On existence of asymptotic solutions to a singular point of an autonomic system possessing a formal solution}, Functional Anal. Appl. \textbf{23}.

\bibitem[KP]{Ref5}
Kozlov V V, Palamodov V P 1982 \emph{On asymptotic solutions of the equations of classical mechanics},
Dokl. Akad. Nauk SSSR \textbf{263} 285--9; English transl., Soviet Math. Dokl. \textbf{25} 335--9.

\bibitem[Lag]{Lagrange}
Lagrange J L 1788 \emph{M\`ecanique analytique}, Veuve Desaint, Paris.

\bibitem[La]{Laloy1}
Laloy M. 1976 \emph{On equilibrium instability for conservative and partially dissipative systems}, Internat. J. Non Linear Mech. {\bf 11} 295--301.

\bibitem[Ly]{Ly}
Lyapunov A M, 1947 \emph{General problem of the stability of motion}, (Kharkov Math. Soc, Kharkov, 1892; French transl.) Ann. of Math. Studies, {\bf 17} Princeton Univ. Press.

\bibitem[LP]{Laloy}
Laloy M, Peiffer K, 1982 \emph{On the instability of equilibrium when the potential function has a non-strict local minimum}, Arch. Rational Mech. Anal. \textbf{78} 213--22.

\bibitem[MN]{Ref6}
Moauro V, Negrini P 1989 \emph{On the inversion of Lagrange-Dirichlet theorem}, Differential Integral Equations \textbf{2} 471--8.

\bibitem[Pa]{Palamodov}
Palamodov V P 1995 \emph{Stability of motion and algebraic geometry}, Transl. Amer. Math. Soc., {\bf 168} 5--20.

\bibitem[Ta]{Ref7}
Taliaferro S D 1990 \emph{Instability of an equilibrium in a potential field}, Arch. Rational Mech. Anal. {\bf 109} 183--94.

\endrefs

\end{document}